\documentclass{article}

\advance\textheight by \topskip
\usepackage{amsmath,amssymb}
\usepackage{color}
\usepackage{graphicx}
\usepackage{amstext,amsthm,amsopn,amsxtra,amsfonts}
\usepackage{setspace}
\usepackage{paralist}
\usepackage[english]{babel}
\usepackage{booktabs}
\usepackage{xurl}

\usepackage{hyperref} 
\hypersetup{colorlinks=true, linkcolor=blue,  anchorcolor=blue,  
citecolor=blue, filecolor=blue, menucolor=blue,  
urlcolor=blue}

\usepackage{color}

\title{Two examples of ungrading in higher education\\from the United States and from Germany}

\author{Christine von Renesse\hspace{0.5pt}\MakeLowercase{$^{\text{1}}$} and\hspace{1.5pt} Sven A.\ Wegner\hspace{0.5pt}\MakeLowercase{$^{\text{2}}$}}

\usepackage[small]{titlesec}
\usepackage[labelfont=bf]{caption}

\begin{document}

\maketitle

\hspace{-1000pt}\footnote{Westfield State University, School of Business, Mathematics, Computing and Sustainability, 577 Western\linebreak\phantom{x}\hspace{10pt}Avenue, Westfield, MA 01086, USA, e-mail: cvonrenesse@westfield.ma.edu.\vspace{1.6pt}}

\hspace{-1000pt}\footnote{Universit\"at Hamburg, Fachbereich Mathematik, Bundesstra\ss{}e 55, 20146 Hamburg, Germany, e-mail:\linebreak\phantom{x}\hspace{10pt}sven.wegner@uni-hamburg.de.}

\vspace{-40pt}

\setcounter{footnote}{0}


\begin{abstract}
In this paper the authors discuss their experiences with ungrading at a small public university in the U.S.\ as well as a large public university in Germany. The courses described are Calculus 1, Mathematics for Liberal Arts, and courses for pre-service secondary teachers of mathematics. We outline and compare our approaches, discuss student performance and feedback and present some interesting patterns relating to gender.

\smallskip

A shortened and revised version of this paper appeared in PRIMUS \textbf{33} (2023), no.~9, 1035--1054, DOI: \href{https://doi.org/10.1080/10511970.2023.2229819}{10.1080/10511970.2023.2229819}.

\smallskip

\end{abstract}

%
%
%
%
%

\section{INTRODUCTION}
The authors met at the virtual joint mathematics conference in 2020 
and started talking about assessment strategies and ungrading in their respective institutions. The different education systems in the United States and Germany put different restrictions on what is possible for a single instructor to achieve and change---which made the comparison even more interesting. 

For us, \emph{ungrading} includes all ways in which the power of evaluation is shared by professor and students. This implies that students evaluate (some of) their work themselves instead of being evaluated by the professor. Dr.\ von Renesse learned about ungrading from Drs.\ D.\ Borkovitz, E.\ Rizzie, I.\ Vasu, C.\ Price, and R.\ DeCoste in the New England Community for Mathematics Inquiry in Teaching (NE-COMMIT \cite{NE-COMMIT}). By now there are several blogs about ungrading \cite{debbie}, \cite{T21a}, \cite{T22}, \cite{C21} and news articles \cite{F19} that we highly recommend reading. Dr.\ Wegner had previously made some experiences with self-evaluation \cite{BW18}. The peer-reviewed articles \cite{Malan21}, \cite{BN86} and \cite{K11} provide additional information about ungrading-related themes. In particular the books ``Ungrading'' \cite{Blum20} and ``Grading the Equity'' \cite{Equity} presented many reasons to change the grading system as well as practical ideas how to do so in various subject areas. To summarize the main idea for change: Students' desire to learn should come from internal motivation instead of wanting to receive a particular grade level. The focus of learning should stay on the content instead of the ``hunt for points.'' Notice that all grading systems are inherently unjust. The COVID19 pandemic has brought to light even more that students need support and how grading policies can harm students' learning and success. While we do not claim that ungrading fixes all the problems with grading, it does break the traditional power dynamics and opens new possibilities. Unfortunately, it also includes new opportunities for bias, as we will present later in this paper.

Dr.\  von Renesse is originally from Germany, but teaches now at a university in the U.S. She implemented ungrading in all her classes in the 2021/22 year, in particular in Calculus 1, Mathematics for Liberal Arts and a content class for future secondary teachers. Because of the pandemic, the classes started online, but moved to in-person after a few weeks.

Dr.\  Wegner teaches at a university in Germany, and implemented ungrading in a content class for future secondary teachers. The course cycle that will be discussed in this article started due to the pandemic online, but went back to in-person teaching in the semester in which the ungrading approach was implemented.

We are both professors of mathematics, not mathematics education, and as such did not prepare a research study about our teaching experience. We did however have questions in mind as we experimented with and discussed ungrading:

\begin{compactitem}

\item[(i)] What effects of ungrading on (a) student learning (b) student grades (c) students life can be observed?

\smallskip

\item[(ii)] How do the differences and similarities between the U.S. and German education system impact ungrading choices?

\smallskip

\item[(iii)] Is there a difference in the ungrading experience based on gender?

\smallskip

\end{compactitem}

To answer our questions, we used qualitative analysis of student journals, minutes of student consultations, our impressions noted during the lessons, and official student feedback of the courses. We also did quantitative statistical analyses of exam and homework results. The classes were relatively small: Dr.\  von Renesse taught 63 students across her courses, Dr.\  Wegner had 40 students in his class.

In the following sections both authors describe and reflect on their ungrading experiences.

\section{UNGRADING AT A SMALL PUBLIC UNIVERSITY IN THE U.S.}\label{SEC-Ungrading-WSU}
Westfield State University has about 4000 mostly undergraduate students and is located in Massachusetts. Many of the students are first-generation college students and work part time or full time to afford college. Dr.\  von Renesse has been teaching at Westfield State University for 15 years. She is teaching all her classes using inquiry, not traditional lecture. 

The following sections describe her ungrading experience in three different courses: Calculus (for Science, Technology, Engineering and Mathematics (STEM) majors), Mathematics for Liberal Arts (MLA), and Mathematical Knowledge for Teaching (MKT). For readers unfamiliar with the U.S.\ system: in most colleges and universities students need to take some mathematics courses for their general education requirement. MLA is one of those classes for majors not in STEM.

\subsection{Ungrading Calculus by Dr.\  von Renesse}\label{SEC-calculus1}
In this section I analyze an implementation of ungrading in a Calculus 1 class in a predominantly undergraduate institution in the U.S.\ in 2022. The 4-credit  class met for 3.5 hours each week and was co-facilitated with an undergraduate teaching assistant. The audience consisted of 18 students from various majors (mathematics, computer science, music, history, biology, psychology, management, Spanish, marketing, chemistry, English, and health sciences). The class was run without traditional lecture: students work in groups on activities that allow them to rediscover the concepts \cite{active-calc}, and whole class discussions are used to introduce vocabulary and notation, as well as sharing our mistakes, insights, connections, and questions. To learn more about teaching with inquiry, \cite{www.artofmathematics.org} is a good starting point. For equity reasons I decided to change the grading system to ``ungrading'' \cite{Equity}. 

During the semester students completed three journals which provided opportunities to reflect on learning mathematics in general: why learning with inquiry is helpful, how our mindset influences our learning, why mistakes are useful, and that there are still inequities in mathematics education in the U.S. The syllabus outlined the following \emph{meta-goals} (not content specific):

\smallskip

Upon successful completion of this course you will:
\begin{compactitem}
    \item Recognize and challenge your own beliefs and feelings about mathematics.
    \item Work well with other students.
    \item Find work partners that have similar speed, curiosity and learning edge.
    \item Gain social competence in negotiating different ways of thinking.
    \item Persevere when it gets difficult and frustrating.
    \item Be ok making mistakes and learning from them.
    \item Be more confident in doing mathematics.
    \item Admit not knowing and ask questions to learn.
    \item Be prepared and ready to try/learn when coming to class.
    \item Realize that mathematics is more about creating and deep thinking than memorizing procedures.
    \item Be curious about mathematics.
    \item Recognize how you learn best and put a plan into action.
    \item Enjoy the challenge of reasoning.
    \item Only accept mathematics that makes sense to you.
    \item Communicate mathematical ideas in writing and oral form.
    \item Reflect on your own experiences as a learner of mathematics.
    \item Explain how there is still injustice in mathematics (STEM) education happening.
\end{compactitem} 

\smallskip

The journals gave me some information about progress toward the meta-goals but not enough to have an influence on the grade level. Students decided how much progress they made on these goals themselves. See \cite{journals} for evidence of effectiveness of meta cognitive assessments.

Homework stories, see description in \cite{repository}, were assigned every 1-2 weeks and graded as ``needs revision'' or ``completed''. Revisions were allowed once a week for each assignment. This meant that the grading load was higher for me since a students could resubmit the same assignment every week until the end of the semester. Fortunately revisions were much easier to grade, since most students just addressed the shortcomings of their prior work and the detailed comments and suggestions helped both parties to focus their attention during the revision process. 

The syllabus contained specifications for content learning goals that had to be achieved for the different grade levels. Many of the content goals were assessed through the homework stories. Notice that the personal and informal writing style of the stories made cheating almost impossible. Two oral conferences provided opportunities for the students to show mastery in more of the learning goals. These meetings could be repeated if necessary. At the end of the semester students used a portfolio to provide evidence for the learning goals they achieved and argued for a final grade. In the final conference the student and I discussed the evidence and agreed on a final grade. The full syllabus can be found in \cite{repository}. 

\begin{figure}
\begin{center}
\includegraphics[width=300pt]{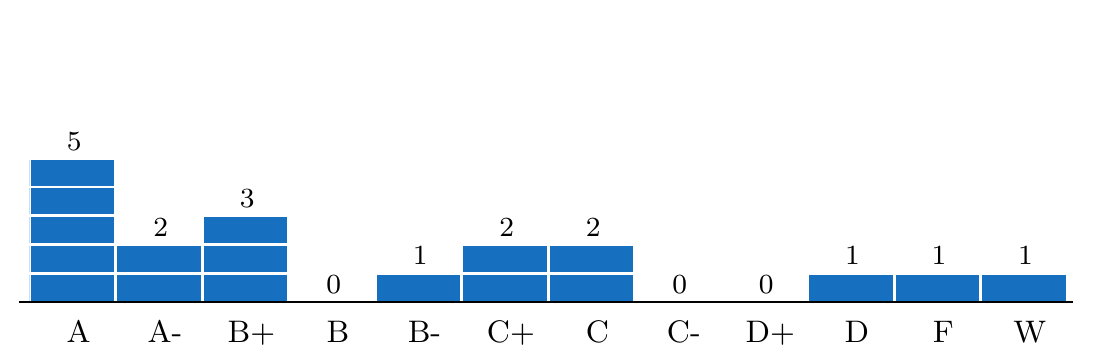}
\end{center}
\caption{Final Grade distribution for Calculus 1}\label{FIG-Calc-Grades}
\end{figure}

One of 18 students withdrew from the class (grade W) and one student stopped doing work and coming to class (grade F). Figure \ref{FIG-Calc-Grades} shows the final grade distribution.

It is interesting to see how the grades that students suggested compared to the grades that I would have chosen. In general the grade choices were remarkably close. Of the six female students, three argued a grade that was lower than the my grade and received a higher grade in the end. Of 12 male students, three argued for a higher grade and agreed on a lower grade during the conference, three argued for a lower grade and agreed on a higher grade during the conference. Only one of the grade difference was more than one letter grade: one of the men wanted a B while I suggested a D+. We agreed on a C in the final conference. 

I noticed a tendency for male students to choose a higher grade than me (24\% of male students, none of the female students), and for female students to choose a lower grade than me (50\% of female students, 25\% of male students). If this pattern happens in general, it could show bias and inequity and it is important for me to be aware of that. Reading the portfolio also gave me insight into the process by which students choose their grades. Notice in the following quote how a female student is working on \emph{confidently} choosing an ``A,'' and how a male student is choosing a ``C.''

\begin{center}
\begin{minipage}{350pt}
``After reflecting on all of these components of the semester overall, I am clearly able to recognize where I will need improvement moving forward in my mathematical career, which this semester I committed to making a longer relationship. The weaknesses I believe I still hold at the end of this semester are being open to communication with those I’m not likely to gravitate towards initially, utilizing online community platforms (aka discussion boards) more to my advantage, practicing more outside of class and keeping practice problems somewhere handy to utilize in homework stories. Past these factors, I am still incredibly proud of myself at the end of the day for my work in this class during the semester and the hardships I’ve learned and lived through these past four months. In conclusion to all of this, I am attempting to \emph{confidently} say that I deserve an A in this class. Overall it has made me work the hardest mentally this semester, bringing me a few tears, laughter, and headaches, and I think I’m taking away the most from this course.'' (female student)
\end{minipage}
\end{center}

\smallskip

\begin{center}
\begin{minipage}{350pt}
``This semester was a genuinely pleasant change compared to last semester. The way you taught the class allowed me to learn a lot more compared to last semester. The thing that worked best for me was definitely the groups, being able to bounce ideas back and forth between peers really makes the 
learning process at least for me a lot easier. For this semester I feel my grade should be a C because I really improved my skills and knowledge this semester, I understand the concepts to a good degree. I feel that I could do all the concepts, I just have trouble with writing them down and making sense of them in words. I feel that I could have done more revisions for my stories and should have been more on top of correcting them and asking you or the discussion board questions about the stories.'' (male student)
\end{minipage}
\end{center}

\subsection{Ungrading Mathematics for Liberal Arts by Dr.\  von Renesse}\label{SEC-MLA}

Mathematics for Liberal Arts (MLA) is a class that satisfies the core requirement of the university. Students come from majors such as nursing, criminal justice, elementary education, psychology, and history. The goals for the class are all meta-goals, the content can be chosen by the professor and/or the students according to interest. Topics included games, Islamic geometry, patterns, and origami, see books at \cite{www.artofmathematics.org}. While the topics are chosen to increase students' engagement and motivation, the emphasis on inductive versus deductive reasoning and informal proofs makes the class mathematically quite challenging. Like the Calculus class described in section \ref{SEC-calculus1}, the class was taught without formal lectures. In addition to the two conferences, students worked on a project of their choice. Instead of exams, students provided a final portfolio to show their success in the class. The homework stories were graded as ``needs revision'' or ``completed'' and revisions were possible on a weekly basis. You can find an example of a homework story from this class with a third revision in our resources folder \cite{repository}.

Students chose their final grade based on the learning goals provided in the syllabus but without any given specifications. I ran the class in prior years with specifications (e.g. for an ``A'' you need to participate at an A level, complete all homework stories, and pass the project, see \cite{specs-blog}) and was curious to see how the full ungrading process would change students' engagement and learning. The only condition was that students had to complete (no further revisions needed) at least half the assignments to pass the class. The full syllabus can be found in \cite{repository}. 

The following student quote showcases that students appreciated the ungrading process:

\begin{center}
\begin{minipage}{350pt}
``The first journal shows where I started with my opinion about math and emphasizes how drastically that’s changed compared to now. I can also recognize that the ungrading process allowed me to focus more on the quality of the work and the energy I put in to get the answer rather than a specific number or letter grade.''
\end{minipage}
\end{center}

\begin{figure}
\begin{center}
\includegraphics[width=300pt]{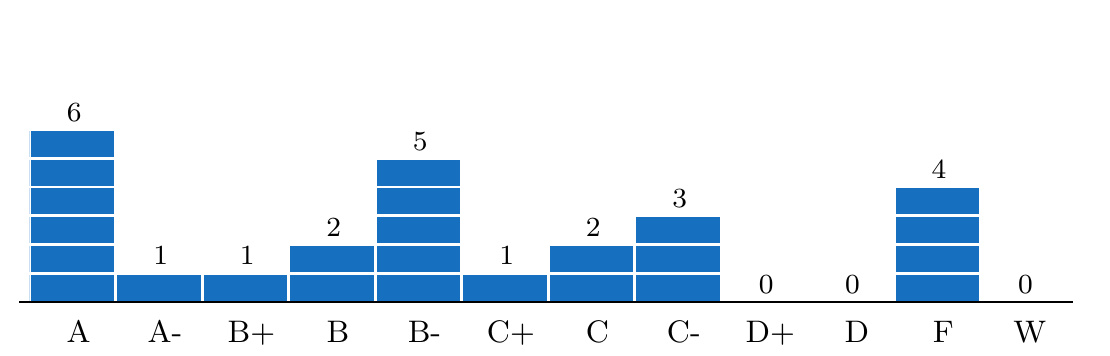}
\end{center}
\caption{Final Grade distribution for Mathematics for Liberal Arts}\label{FIG-MLA-Grades}
\end{figure}

The final grade distribution of the MLA class is shown in Figure \ref{FIG-MLA-Grades}. Notice that it became clear from the conversations during conferences that C is the worst grade a student would choose. Students told me that they considered anything worse than a C a ``fail''.

Again, the gender pattern was interesting to observe: 12 of the 25 students are female and 13 are male. Three female students stopped attending class, one female student didn't do most of the work needed for the class. Of the remaining eight female students, three argued for a grade that was lower than my grade and received a higher grade in the end. The other five students chose the same grade as I did. Of the 13 male students, nine argued for a grade that was higher than the my grade and received a lower grade than they had originally chosen. Two male students had a grade difference of two between their choice and my choice. As in the Calculus class in section \ref{SEC-calculus1}, I notice a tendency for male students to choose a higher grade, and for female students to choose a lower grade compared to my grade suggestion. 


\subsection{Comparing Calculus and Mathematics for Liberal Arts by Dr.\  von Renesse}\label{SEC_CompareCalcMLA}

In both classes, about half the students chose to revise each assignment (on average). This is interesting since Calculus had the added incentive of given specifications for each grade level and Calculus is a prerequisite for STEM majors for other courses, while MLA only had a minimal threshold for passing the class. I would have guessed that the Calculus students would do more revisions than the MLA students.

It seemed to me that students who revised a story, did so to either learn or to guarantee getting an A. Students did not revise to avoid failing or getting a bad grade so the “punishment” effect of a minimal threshold did not do much in comparison to the fall semester (when I ungraded without a minimal threshold). Maybe this was the case because students overall were much more stressed and overworked in comparison to pre-covid years? 

In both classes, there were about 60 revisions of homework stories total (counted for each student and each revision separately). Both classes had about 20 active students so are comparable. For six homework stories, I would get at most 120 items to grade without revisions, now, with revisions, there are 60 more items to grade---a 50\% increase. I used to assign a story for each course every week but since I started the revision process, I reduced the assignment load by 50\% to keep a doable workload for myself and the students. 

There were a few stories for which students received feedback in class. These real time revision suggestions were not recorded in the grade book. Overall students revised earlier stories more often than later stories probably because the later stories had less time for revision.

Notice that in MLA women did more revisions than men. 38 revisions were done by women and only 23 by men. In Calculus the number of revisions by gender was about the same. This pattern may have been influenced by student majors in the different classes.

For the final 10 minute conference in MLA, I started out by showing the students our list of meta-goals and asking them to reflect on their progress. After that I would tell them what I thought about their portfolio, and we talked about their final grade choice and justification. In Calculus, the 20 minute conference included an oral exam period: For the first 10 minutes I asked them questions about the class materials. The second half was spent like the MLA conferences. 

I noticed that overall male students seem to present themselves with more confidence even if they had less successful work to show for than their female counterparts. I often helped female students see work they had accomplished (although not always perfectly), and asked them what they thought they should have done for a higher grade. Often the answer was: I should have struggled less, or I should have done work on time more often, or other judgements that in my opinion should not influence the final grade too much. On the other hand I presented some of the male students with non-finished work, and asked them afterwards to reassess their grade. When we could not agree on anything, an “incomplete” was an option to finish more work for the higher grade level. 

When students did not have time to revise a story again but felt they fully understood the topic, I often counted a topic or story as completed even without evidence. This influenced some of the final conference grade decisions. I decided to trust the students, thinking that their assessment of what they understand at that point was better than mine, since I had no evidence and they had some.

\subsection{Ungrading Mathematical Knowledge for Teaching by Dr.\  von Renesse}\label{SEC-MKT}

At Westfield State University, future secondary teachers take all classes of a regular mathematics major. They are taking additional classes in education and psychology. In addition, they are taking several classes specific to the teaching of mathematics, like Mathematical Knowledge for Teaching (MKT). The focus is on understanding the content of middle and high school more deeply, including proofs, counter examples, common mistakes and misconceptions. In addition, the class focuses on how to teach mathematics: through modeling in class, class visits of freshmen college classes and reading assignments. Students prove, for example, the quadratic equation using algebra tiles, the equation for the volume of a cone using Cavalieri's principle, and Pick's theorem.

Like the MLA class described in section \ref{SEC-MLA}, the class was taught without formal lectures. In addition to the two conferences, students visited other classes and reflected on the teaching they observed. Instead of exams, students provided a final portfolio to show their success in the class. The homework stories were graded as ``needs revision'' or ``completed'' and revisions were possible on a weekly basis. Students chose their final grade based on the learning goals provided in the syllabus but without any given specifications. The only condition was that students had to complete (no further revisions needed) at least half the assignments to pass the class. The full syllabus can be found in \cite{repository}. 

\begin{figure}
\begin{center}
\includegraphics[width=300pt]{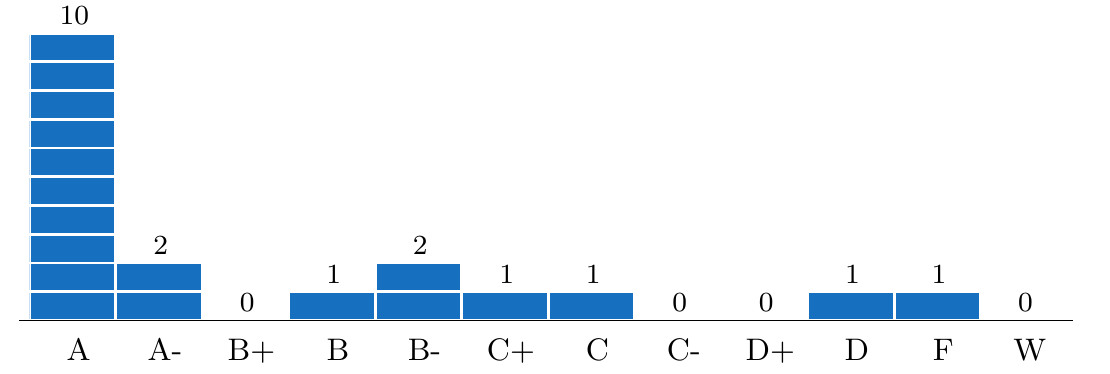}
\end{center}
\caption{Final Grade distribution for Mathematical Knowledge for Teaching.}\label{FIG-MKT1}
\end{figure}

The final grade distribution for the MKT class is shown in Figure \ref{FIG-MKT1}. Notice that students in this class were very invested in learning the topics, which helps to explain why so many students received high grades.

The class consisted of 10 female and 9 male students. Three women graded themselves higher than I would have, one graded herself lower. In total two of the women received a grade higher than what they chose, one received a lower grade. All grade differences were less than one letter grade. 
Two of the men argued for an incomplete and received a grade a few weeks after classes ended. One man stopped attending class and received an F. Three men graded themselves higher than I would have, while for one man I chose a higher grade. In total three men received a higher grade than what they chose, one received a lower grade. All grade differences were less than one letter grade. 

The gender pattern of the other classes can not be observed here. However, seven women received an A, and two women received an A- so the female students were stronger in the class than the male students. 

For this class, 12 students revised each story on average, so about 60\% of students. Since one student was basically not present, it is more accurate to say that 66\% of active students revised each story at least once. Of the 91 total revisions, 62 were done by women, and 29 by men. So while the grade changes do not show a difference by gender, there is a difference in revisions done. 	

The following quote from a student portfolio, showcases how students realize why ungrading is helpful for learning:

\begin{center}
\begin{minipage}{350pt}
``Throughout this course I honestly struggled more than I expected to. I felt like it would be simple but once this new format of grading actually proved to be a challenge for me because it really focused on understanding every part of the assignment rather than getting a certain percentage of it correct and not worrying about the rest. When it comes to math you will not get very far totally understanding like 70\% and totally ignoring the other 30\%. You truly need to have a good grasp on topics to have a solid foundation to build from and the grading system really reflects that. In my first journal I even talked about how I do not personally agree with the ungrading system. I guess just reading it did not show me truly the benefit of it. [...] As far as grades go I think I deserve an A in this course because I feel like I not only meet the expectations of the course but I pushed myself to truly understand more than just a grade. I think I have really understood what ungrading is about and realize the grade shouldn't be what is important but rather actually taking the time to learn the material.''
\end{minipage}
\end{center}

In terms of workload, five stories would be 95 items to grade in total with no revisions. The revisions process added 90 more items to grade, so the grading load for me roughly doubled.

\subsection{Reflecting on the Ungrading Experience by Dr.\ von Renesse}\label{SEC-reflection-CvR}

Reading the book ``ungrading'' \cite{Blum20} convinced me to try something similar in my mathematics courses. The pandemic had changed so many things already and showed more clearly the many ways our grading system is not equitable. In my regional COMMIT group \cite{NE-COMMIT} we read the book ``grading with equity'' \cite{Equity} which added more perspectives. I had used specification grading in all my courses for the three years prior to fall 2021. Part of the specifications were meta-objectives and participation rubrics. The ``grading for equity'' book convinced me to not grade ``behavior'' anymore, and the ungrading perspective allowed me to have students choose goals and evaluate their progress without me judging them.

For example, a quiet student may want to improve their participation in whole class discussions and can then reflect on the progress they made. Or they can choose not to do that and instead change their attitude towards making mistakes in mathematics and  feeling confident on their journey. Either way, a free writing assignment in class, several journals, the conferences, and the final portfolio give the student ample opportunity to choose reasonable personal goals and work toward them. It was a pleasure to observe how students changed over the semester and how they reflected on their growth. Of course there were a few students in each class who didn't attend class enough, or didn't do the assigned work---just like in any of my classes over the years. But students mostly chose to come to class, and do the assigned work. Since I do not have strict deadlines, some students procrastinated and in hindsight wished I had given them strict deadlines. But most students were grateful that my class didn't add as much extra pressure as other classes, and found time to eventually do the work. Attendance in my classes was lower than attendance in classes were absences were punished or graded, of course. There were days when I was frustrated by students not coming or showing up late. But then I heard students share what got in the way (physical or mental illness, a sick family member, a sibling that needs someone to watch them, wedding preparations, no sleep the night before because of an exam,...) and I was so glad that they did not force themselves to come to class. 

Since the pandemic all classwork done by groups of students during class is posted on google docs.  Absent students can attend class synchronously via zoom or participate asynchronously by working on the google docs by themselves. So not being in class did not necessarily mean that students did not engage in the material.

The student population at my university has changed over the last years due to the pandemic and due to a student-shortage in New England. Overall students at my school are less prepared in mathematics, and for the college experience in general. Ungrading a class allows me to focus on supporting students in their learning instead of playing the game of grades and looking for ways to punish/reward them. It was a pleasure to have productive conversations about content, study habits, and life in general---without the distractions of grades.

I had to learn to trust my students, instead of just expecting them to trust me. At times I find this difficult to do, but open conversations and being vulnerable myself, usually leads to the students showing up more and differently. Ungrading allows for connection and trust. In my opinion, learning should be engaging and fun (which doesn't mean easy!)---not regulated by constant fear and stress. 

In classes that are pre-requisites for other classes, like Calculus 1, I do like the mixture of some specifications, some ungrading for content objectives and complete ungrading for meta objectives. I can not, in good conscience, allow totally unprepared students to move on to the next class in a sequence like Calculus.

For my classes next semester, I am planning to continue ungrading in some way. The revision process added extra work when I introduced it years ago, but assigning fewer stories helped with that. The conferences take a huge amount of time though. With 80 students last semester in four classes I had a difficult time arranging the meetings and spent many evenings on zoom. I do not think the one-on-one conferences would be feasible for me with more than 80 students.

I need to point out that all my classes are small and that any evidence I provided is limited in scope, especially concerning the gender differences. I hope though that my data will start an interesting and challenging conversation about learning and teaching which will in turn lead to more research and ideas. 

I want to finish my reflection with a quote of a student from one of my graduate classes for teachers. While this class was not discussed in this paper, the reflection shows the power of ungrading---so that I chose to include it here. My graduate students have much more experience in reflective writing than my undergraduate students, but I heard the same sentiment from many of my undergraduate students.

\begin{center}
\begin{minipage}{350pt}
``Taking everything that I have stated so far [in my portfolio], and reflecting on the effort that I put into this class, I would give myself an A for this class. My initial reaction when I learned that I would be evaluating myself for a final grade was one of excitement. Then as the class was progressing I started to get nervous because I wasn't sure if I was doing “enough” for the grade that I wanted. I was doing all of the readings, participating on Discord and trying out the new concepts in my classroom, but I still wasn't sure it was enough. However, after reading the “homework” chapters in both of our books and reflecting on the whole purpose of “homework/practice/learning in general being” for the student and not the teacher, I kind of changed my perspective on things. I started to reflect on what I was truly getting out of this class and whether or not the effort I was putting in was giving me the results THAT I WANTED, not necessarily what was expected of me by anyone else. A funny thing happened when I did that... the pressure of needing to excel and prove myself went away and I became more willing and vulnerable in what I was willing to try and risk in the classroom. Then, as I started putting together all of my portfolio pieces, I realized just how much I got out of the class and all of the effort that I put in along the way. I also started to think about how much learning I am still planning to do in terms of making my classroom a “thinking classroom” moving forward and all of the steps that I have started to take to get there and the effort and groundwork that I have already put into it. It was at that point that I felt I was comfortable advocating for an A in the class, not just because I wanted it, but because I earned it.''
\end{minipage}
\end{center}

\section{UNGRADING AT A LARGE PUBLIC UNIVERSITY IN GERMANY}\label{SEC-Sven}

The University of Hamburg (UHH) has about 45,000 students and offers a wide range of degrees in various subject areas. It is one of the eleven ``excellence universities'' in Germany. Dr.\ Wegner is since 2020 a Tenured Research Fellow at UHH with teaching duties in the programs for pre-service high school teachers. He came to Hamburg with more than 10 years experience in graduate and undergraduate teaching of which he spent 2.5 years at a post-92 university in the UK and one year at a start-up university in Kazakhstan which was running a U.S.-style liberal arts model.

In this section Dr.\  Wegner analyzes an implementation of ungrading in the course cycle Mathematics 1--4 taught between 2020 and 2022. The four courses are the mandatory foundational mathematics courses for pre-service high school teachers and usually are taken in the first two years. In Germany, high school teachers are required to teach two subjects (e.g.~mathematics and physics, mathematics and biology, mathematics and English, etc). Their university education thus comprises courses on these two subject areas, courses on general pedagogy and courses on subject specific pedagogy. Traditionally, pre-service high school teachers of mathematics enrolled within the courses Analysis 1--2 and Linear Algebra 1--2 together with regular mathematics majors. A change in the study regulations lead to the new courses Mathematics 1--4 taught only for pre-service teachers. The incentive behind this change was to reduce drop-out rates and increase the quality of mathematics teacher education. The 2020 cohort that was part of Dr.\ Wegner's ungrading approach was the first who enrolled in the new courses.

Dr.\ Wegner was the main person responsible for tutorial organization of Mathematics 1, 2, and 4 which were taught by the UHH professor Dr.\ Thomas Schmidt. The course Mathematics 3 was taught by Dr.\ Wegner himself. In all four courses there was support by several student assistants.

\subsection{Ungrading Mathematics 3 for Pre-Service Teachers by Dr.\ Wegner}
\subsubsection{Course Background}\label{SEC-Sven-1} 


\smallskip

All four courses in the cycle Mathematics 1--4 have a value of 9 ECTS credits (roughly corresponding to 4--5 credits in the U.S.) and consist of two 90-minutes lectures (given by the professor in charge), one 90-minutes tutorial (where weekly homework is presented by students themselves or by a teaching assistant) and one 90-minutes workshop (where students work on tasks in small groups (supervised by two teaching assistants). For lectures, all students are taught together while for tutorials and workshops they are split up into groups of size less than 25 students each. In view of its volume, the cycle naturally covers a wide range of topics. From the U.S.\ perspective, and only as a rough point of reference, one can think of the courses Mathematics 1--4 as an amalgamation of upper-level courses like Foundations of Mathematics, Introduction to Proofs, Abstract Algebra, Advanced Linear Algebra, Real Analysis I and Real Analysis II into one giant course ranging over sophomore and junior year.\footnote{German universities do not have a freshmen year. In mathematics that means that students go from high school calculus and algebra directly to proof-based courses. A Bachelor's degree in Germany requires three years of study.} Notice that this course, while also serving future secondary teachers, is very different in content from the class that Dr.\ von Renesse taught, see section \ref{SEC-MKT}.  Dr.\ von Renesse's class covers content from the high school curriculum, while Dr.\ Wegner's class covers the typical content for a mathematics major at the university level. A detailed description of module contents and study aims at UHH can be found in the official documents \cite{FSB, UHH-Mod}.

\smallskip

The core of tertiary mathematics education in Germany are traditional lectures which go along with weekly homework sheets that are then discussed in the tutorials. The homework sheets usually contain computational tasks as well as proof-based tasks that focus on understanding. For the courses Mathematics 1--4 the concept of a weekly ``workshop'' was added. Here students work in small groups on problems that (a) process the topics of the lecture, (b) prepare the homework tasks and (c) connect the topics of the course to school mathematics. The idea behind this is to dovetail the subject courses in mathematics with the subject pedagogy courses and create a more coherent education of pre-service teachers.
Although this approach is not exactly a bottom-up approach it has strong elements of 
discovery/inquiry since students have to come up with their own conjectures, refute arguments, find mistakes in calculations and experiment with definitions, see  \cite{PW20}. 

\smallskip

The course cycle started with a cohort of approximately 55 students of which approximately 35 finished the cycle with us. The remaining students either retook parts of the cycle, paused, or dropped out of the program. Due to the COVID-19 pandemic, Mathematics 1 and 2 were taught remotely. Mathematics 3 was taught in a hybrid mode, that is, the lectures were given in the classroom but a recording was made available afterwards. Every week one tutorial and one workshop were taught remotely, while the other tutorials and workshops were done in-person. In Mathematics 4, lecture recordings were still made available, but all tutorials and workshops were done in-person. Throughout the whole cycle, students had to submit the weekly homework sheets in teams of two (and in one exceptional case of three) students. Concerning assessment we had the following systems: In Mathematics 1 and 2 students received points on each assignment and were required to score 40\% or above on the homework in order to be eligible to participate in the final exam. In the case of Mathematics 1 this was, due to the pandemic, a pass/fail oral group exam while in Mathematics 2 it was an in-person written exam in which 35\% of the points were sufficient to pass.\footnote{In German tertiary mathematics education it is a very common setup to have in-course-assessment as an entrance requirement for the end-course-assessment which in turn has the form of an exam which then solely determines the (letter) grade. Indeed, study regulations for our courses Mathematics 1--4 prescribe that the end-course-assessment has to be a written exam. The oral group exam in Mathematics 1 was an exception granted by the exam committee due to the pandemic.} After having learned from Dr.\ von Renesse in the spring of 2020 about ``ungrading,'' I decided to give the concept in Mathematics 3 a try and I altered the paradigm as follows: There were still weekly homework sheets to be submitted by students in teams of two and they were checked by a teaching assistant but the grading was dropped. More precisely, students received mathematical feedback on errors in their submission, but no points were assigned by the teaching assistants. Instead, students themselves were required to assign for each problem either ``ok'' or ``needs work'' and, additionally, to reflect and comment with a few lines on their learning progress. The requirement to score 40\% on homework for participating in the final was replaced by individual end-of-term consultations of me and submission team, in which exam participation was discussed. In Mathematics 4 Dr.\ Schmidt, who was then again the professor in charge, returned half way and re-introduced non-numerical grades ``exemplary,'' ``correct,'' ``half-correct'' and ``wrong'' given by the teaching assistants on each assignment. Indeed, he asked the students in the first lesson what they would prefer and in view of their responses he agreed on this compromise. He kept however the ``self-evaluations'' and did not re-impose the 40\% threshold but instead communicated to students that consultations would be undertaken only in cases were professor and teaching assistants determined that this was necessary. Indeed, all but two students (who stopped submitting in week 3) were admitted to the exam without consultation.

\smallskip

The next section will focus on the courses Mathematics 2 and 3. The aim is to analyze the impact of my version of ungrading (self-evaluation of homework + consultation + teacher-graded exam) in Mathematics 3 and then compare it to the classically taught Mathematics 2 (teacher-graded homework + 40\% threshold + teacher-graded exam).

\subsubsection{Data Acquisition and Analysis}\label{SEC-Sven-2}

As mentioned in section \ref{SEC-Sven-1} I introduced my version of ungrading at the beginning of Mathematics 3. In the first lesson I explained the concept, namely that (a) homework will still be checked and feedback will be made available but no points will be given, and that (b) participants are supposed to assign themselves either ``ok'' or ``needs work'' for each assignment and to write a comment on their learning progress. The latter process was labeled ``self-evaluation.'' I emphasized several times throughout the course that the idea of  the latter is to re-engage with the topic of the assignment \emph{after} mathematical feedback was received and solutions were discussed during the tutorials in order to foster understanding in the long run.\footnote{In the actual course, I asked students to write P$\,\equiv\,$PASS for a task that they mastered and F$\,\equiv\,$FAIL for a task where they feel that it will need more attention, before the exam at the latest. Retrospectively, my impression is that the labels ``ok'' and ``needs work'' would have been much better, in particular as F comes with the connotation of an (irreversibly) failed assignment, whereas the term ``needs work'' underlines an opportunity for improvement. See also \cite{T21}. In the remainder of this article we stick to ``ok'' and ``needs work'' in order to propagate these notions among our readership.} Moreover, the self-evaluations were supposed to help during the exam preparation in that students can then (a) use their own comments to refresh their memories more quickly when revisiting tasks they had mastered before and (b) invest the saved time to focus on tasks labeled ``needs work.'' While some students followed this advice closely, many others wrote rather generic comments, apologized for not submitting or submitting something incomplete, or wrote no comment at all. I classified the comments in four categories, namely ``insightful,'' ``generic,'' ``apologetic'' and ``no comment.'' Here, insightful comments are comments that contain a mathematical insight specifically related to the task. Generic comments are comments of general mathematical nature or no mathematical reference at all. Purely descriptive comments also fall in this category. Apologetic comments are comments in which non-mathematical reasons are given why a task could not be attempted or completed. With ``no comment'' we label those responses containing only ``ok'' or ``needs work'' and nothing further but not tasks with no response at all. You can find a list of a few examples of comments in the first three categories below:
\vspace{3pt}

\begin{center}
Generic:

\vspace{3pt}

\textquotedblleft{}Imprecise mathematical arguments\\(unclear what we wanted to say), formal errors.\textquotedblright{}

\vspace{3pt}

\textquotedblleft{}We understood the task and did answer it correctly.\textquotedblright{}

\vspace{3pt}

\textquotedblleft{}Problem with understanding, therefore wrong result.\textquotedblright{}

\vspace{3pt}
 
\textquotedblleft{}I understood this task well after the tutorial,\\but I intend to go over it again before the exam.\textquotedblright{}

\vspace{9pt}

Insightful:

\vspace{3pt}

\textquotedblleft{}From the tutorial it is now clear how to use the definition of convergence\\of sequences and by defining $N:=\operatorname{max}\{N_f,N_g\}$ [\dots] can conclude\\that the maximum [of the functions $f$ and $g$] is continuous.\textquotedblright{}


\vspace{5pt}

\textquotedblleft{}We overlooked that gaps in the domain result\\in a function that is not anymore defined on $[a,b]$.\textquotedblright{}



\vspace{9pt}

Apologetic:

\vspace{3pt}

\textquotedblleft{}We had an important presentation this week,\\therefore we could do only one [homework] task.\textquotedblright{}

\end{center}

\vspace{5pt}

In Table \ref{TAB-2}, you can find the percentages of comments in each category throughout 55 homework tasks given in the course Mathematics 3. It also shows the percentages if we differentiate between the labels ``ok'' and ``needs work'' given by the students. For the labels I need to point out that some students made up intermediate labels by writing e.g. ``something in between ok and needs work.'' In such cases I counted the task as ``ok.'' For tasks with several parts, some students assigned a label to each part. Here I used a majority vote and broke ties in favor of ``ok'' to determine how to count the task. There were 19 submission teams of which one submitted only two homework sheets and no self-evaluations. This team was therefore be excluded from the analysis in Table \ref{TAB-2}. One submission team did the self-evaluation not as a team but for each of the two members separately.  The other teams, consisting of two and in one case of three students, did the self-evaluation team-wise. The table contains data of $55\times(17+2)=1045$ self-evaluations given by $17\times 2 + 1\times 3=37$ students. The reason that the percentages do not add up is due to rounding errors.

\vspace{5pt}
\begin{table}

\begin{center}
\resizebox{\textwidth}{!}{\begin{tabular}{cccccccccccccc}
\toprule
    \multicolumn{2}{c}{insightful}    && \multicolumn{2}{c}{generic} && \multicolumn{2}{c}{apologetic} & &\multicolumn{2}{c}{no comment}& & \multicolumn{2}{c}{no response}  \\
   \cmidrule(lr){1-14}
       \multicolumn{2}{c}{total} &	&\multicolumn{2}{c}{total}  & & \multicolumn{2}{c}{total}   & &\multicolumn{2}{c}{total}   & & \multicolumn{2}{c}{total}	\\
    \multicolumn{2}{c}{10.3\%} 	&& \multicolumn{2}{c}{55.7\%}  & & \multicolumn{2}{c}{16.7\%}  & & \multicolumn{2}{c}{4.2\%} &&  \multicolumn{2}{c}{13.2\%} 	\\
     \small ok &\small n/w  && \small ok &\small n/w  && \small ok &\small n/w && \small ok &\small n/w  \\
   \phantom{0}7.8\% 	&\phantom{0}2.5\% & & 37.3\%  & 18.4\% && 15.4\% &\phantom{0}1.3\% && \phantom{0}1.1\% & \phantom{0}3.0\% 	\\
\bottomrule
\end{tabular}}

\end{center}

\caption{Labels and comments by students in Mathematics 3.}\label{TAB-2}

\end{table}

\vspace{5pt}

At first, the data looks disappointing since there are rather few insightful comments. I want to point out however, that it cannot be expected that a majority of comments are insightful: If a student has completed a task without mistakes then my system does not require to write anything. If a student did not attempt a task due to time constraints then there is likely no epiphany  during the tutorial when someone else presents the solution. On the other hand it may happen that a student worked hard on a task, did not complete it, had a major insight during the tutorial but simply did not write anything about that down.

\smallskip

In view of the above I will now focus on four cases in which I will, in connection with a specific topic of the course, analyze (1) a homework assignment, (2) the corresponding self-evaluation, (3) the discussion during the consultation and (4) a similar or related assignment in the exam.

\smallskip

\textbf{Case Study I on uniform continuity:} The concept of uniform continuity was defined right after continuity. Examples were given in class and the theorem stating that continuous maps on compact intervals are uniformly continuous was proved. On the corresponding homework sheet we had the following two assignments:\vspace{3pt}

\begin{compactitem}

\item[1.] Let $f\colon(a,b)\subset\mathbb{R}\rightarrow\mathbb{R}$ be continuous with $a<b$ real. Show that $f$ is uniformly continuous if and only if there is a continuous extension of $f$ to $[a,b]$.\vspace{3pt}

\item[2.] Show that $f\colon(0,1]\rightarrow\mathbb{R}$, $f(x)=\sin(1/x)$ is not uniformly continuous.\vspace{3pt}

\end{compactitem}
The idea here was that students can apply the statement in 1.\ to the concrete situation given in 2. We received  submissions for 1.\ from 11 teams of which 9 were correct. For 2.\ we received 11 answers which all were correct. However, all those answers checked by hand that $f$ is not uniformly continuous neglecting the statement in 1. Although the tutorials covered this and it was pointed out, that it is not necessary to disprove ``by hand'' that $f$ in 2.\ is not uniformly continuous, no self-evaluation contained an insightful comment. In the consultations I discussed the above again with all groups. Some expressed surprise and had apparently not memorized the connection between the proof-based task and how it can be useful in concrete situation. I pointed out that students should anticipate a question about uniform continuity in the exam. There, I asked:\vspace{3pt}

 \begin{compactitem}

\item[3.] Determine where the function $f\colon(0,1]\rightarrow\mathbb{R}$, $f(x)=\sin(1/x)$ is continuous and check also if it is uniformly continuous.\vspace{3pt}

\end{compactitem}

The second part of this question was attempted by all 26 exam participants and done successfully by 12 of them, where eight employed the criterion from 1. Two more used the criterion but did not complete the task correctly. This indicates that the consultations might have had a positive impact on students' learning progress concerning this topic.

\medskip

\textbf{Case Study II on uniform convergence:} The concept of uniform convergence of sequences of functions had been covered already in Mathematics 2 in the context of power series. In Mathematics 3 we repeated the definition, gave a few examples and then proved, among other statements, that if a sequence of continuously differentiable functions converges pointwise and the sequence of its derivatives converges uniformly then the pointwise limit function is continuously differentiable and limit and derivative commute. We had the following two questions combined as one task on the corresponding homework sheet:\vspace{4pt}

\begin{compactitem}
\item[1.] Show that $(f_n)_{n\in\mathbb{N}}$ defined via $f_n\colon\mathbb{R}\rightarrow\mathbb{R}$, $f_n(x)=\sqrt{1/n+x^2}$ converges uniformly to the absolute value $|\cdot|\colon\mathbb{R}\rightarrow\mathbb{R}$. Does the sequence $(f_n')_{n\in\mathbb{N}}$ of derivatives converge uniformly, too?
\end{compactitem}
\vspace{5pt}

This assigment was attempted by 16 submission teams, of which 5 submitted a correct solution. The other 11 teams either used the supremum norm incorrectly (they interchanged limit and supremum \emph{to show} uniform continuity) or did only consider the pointwise limit. In the tutorials this was discussed and indeed here we got the following insightful comments in the self-evaluation:

\vspace{5pt}

\begin{center}
\textquotedblleft{}We interchanged limit and supremum which is in general\\not valid and has to be proved to be true in special cases.\textquotedblright{}

\vspace{5pt}

\textquotedblleft{}For the first question one can argue via pointwise\\ convergence of $f_n$ to $|\cdot|$: Since  $|\cdot|$ is not differentiable\\in $0$, $f_n'$ cannot converge uniformly, cf.~Theorem 10.12.\textquotedblright{}


\vspace{5pt}

\textquotedblleft{}In the first task we had to show the uniform convergence of\\the derivatives. We thought, that the derivatives of the function\\needs to converge to the derivative of the absolute value.\textquotedblright{}


\end{center}

\vspace{5pt}

We discussed this task and the use of the supremum norm in the consultations and I emphasized that in the exam there will be a question on uniform convergence as well. The exam contained the following task:
\vspace{5pt}

\begin{compactitem}
\item[2.] Show that a sequence $(f_n)_{n\in\mathbb{N}}$ of function $f_n\colon D\subset\mathbb{R}\rightarrow\mathbb{R}$ converges uniformly to $f\colon D\rightarrow\mathbb{R}$ if and only if $\displaystyle\lim_{n\rightarrow\infty}\|f_n-f\|_{\infty}=0$ holds.
\vspace{5pt}
\end{compactitem}
The above fact was mentioned in the lectures (and the equivalence had been used in other theorems and in the homework that we explained above) but the proof had not been given. In the exam there were, out of 26 submissions, two almost perfect solutions, two half correct solutions and 24 solutions with nothing or almost nothing correct. In particular, from the students that gave the three comments above, only one delivered a half-correct solution. The two almost perfect solutions stem from students who did not submit an insightful comment. Of those two, one performed already well in 1., the other did not. This case analysis does not indicate an uptake in learning outcomes caused by my ungrading approach. I would like to mention however, that from my experience, uniform continuity is a topic that many students find difficult and where scores are also often very low in a traditional assessment setup.

\medskip

\textbf{Case Study III on Taylorexpansion:} I presented the usual theory of one-dimensional Taylor expansion in the lecture, including the Lagrange remainder. As part of the homework I assigned  the following task:\vspace{5pt}

\begin{compactitem}
\item[1.] Compute the Taylor polynomial $\operatorname{T}_{3}[x^2+3x+4,10]$.
\vspace{5pt}
\end{compactitem}

Indeed, 18 teams did this correctly up to computational errors, but all computed it by hand. In the tutorials it was explained how this task can be completed without any computation. However, no student made a comment on that in the self-evaluation. In the consultation we did not discuss this explicitly again, but I pointed out that for the exam one should know how to compute Taylor polynomials by hand and by using ``tricks'' that were discussed in the homework. The corresponding exam task was: \vspace{5pt}

\begin{compactitem}
\item[2.] Compute the Taylor polynomial $\operatorname{T}_{10}[x^5-3x^4+7x-1,1]$.\vspace{5pt}
\end{compactitem}

Indeed, here all 26 exam participants attempted it, 15 completed it correctly of which five did it without computing. Nine students had the correct computational approach but made computational errors. In an informal discussion with some students during Mathematics 4 (when the same question arose in the context of multidimensional Taylor) several students said that in the Mathematics 3 exam they were unsure if it ``was allowed'' to argue $\operatorname{T}_{10}[f,1]=f$ without doing the computation. Here in particular, I believe that a suitable comment in their self-evaluations would have helped those students much---if they would have written some.

\medskip

\medskip

The case studies above illustrate how I used a self-evaluation technique to ungrade the course Mathematics 3. Case Study I suggests that positive effects on the learning progress can happen (a) through meaningful insights that the students understand more deeply by writing about them in their own words and (b) via concrete feedback given during the consultations that points out specific mistakes and ways how to improve. On the other hand, Case Study II shows that even after insightful comments have been made by students and topics had been revisited during consultations, it may still happen that exam performance on a corresponding task is very low. Case Study III illustrates that an increase in performance from homework to exam can happen also without insightful comments or concrete pointers during  consultations. You can also see that students did not take advantage of the system to its full potential (in particular in Case Study III it would have been very easy to write a helpful comment, and by adhering to it in the exam, to save time and avoid point deduction).

\smallskip

As I mentioned in section \ref{SEC-Sven-1}, the admission threshold of scoring more than 40\% on teacher-graded homework in Mathematics 2 was in Mathematics 3  replaced by the consultations. For the latter every submission team had to self-register for their consultation slot. The one team that had submitted only occasionally did not register for a consultation. The other 18 submission teams participated in the consultation. Prior to the consultations students were informed that they had to (a) present one exercise that they found particularly hard at first but then clear and easy in the end and (b) present one exercise or topic that they saw already in high school but where this course significantly extended their understanding. I used this as an ice breaker to begin the discussion. After that I asked if there were any assignments that students wanted to talk about again. A few groups brought up topics that we then discussed and where we could clarify misunderstandings. Finally, we worked through an exam checklist, which had been handed out to students in advance, and focused on six list items relating to six homework questions. For each team I had, prior to the consultation, looked again at their corresponding submissions and I tried to give concrete pointers for the exam preparation. Examples are explained in the Case Studies I--III. The consultations took approx.\ 40 minutes each. For 17 teams it was a productive meeting that had the character of a coaching session. There was only one team where my consultation's preparation had shown that this team had submitted only less than a third of the assignments and of these again a third done only partially. Since I saw no chance for this group to pass the exam, I decided not to admit the students and communicated this during the consultation.

\smallskip

Next, I want to compare the homework labels and the exam results. For that, I counted in Mathematics 3 the number of ``ok'' labels and the number of ``needs work'' labels for each submission team. For Mathematics 2 I used the teacher-assigned score. To compare this with the individually written exams, all tables, histograms and plots below show, with respect to homework, for each student the team's results. In case that a student did the first attempt of an exam, failed and did then the second attempt as well, the exam score is computed as the maximum of the two scores. I want to point out again that Mathematics 2 and 3 were taught by different professors in charge. Both courses were passed by passing the exam. In  Mathematics 2 the pass threshold was 35\% and in Mathematics 3 it was 50\%. This does not mean that the level of the courses is different but simply reflects different exam designs.

\smallskip

First, I compare admission and pass rates of Mathematics 2 and 3: In Mathematics 2, 90\% of students were admitted to the exam. Of those, 86\% participated in the exam and of those 86\% passed the exam. The overall pass rate of the course was thus 68\%. In Mathematics 3, 92\% of students were, as a result of the consultation, admitted to the exam, of these 88\% participated and 93\% of those passed. This gives an overall pass rate of 76\%. In particular, all students who gave themselves ``ok'' on more than 40\% of the tasks and participated, passed the exam. 

\smallskip

\begin{figure}
\begin{center}
\includegraphics[width=190pt]{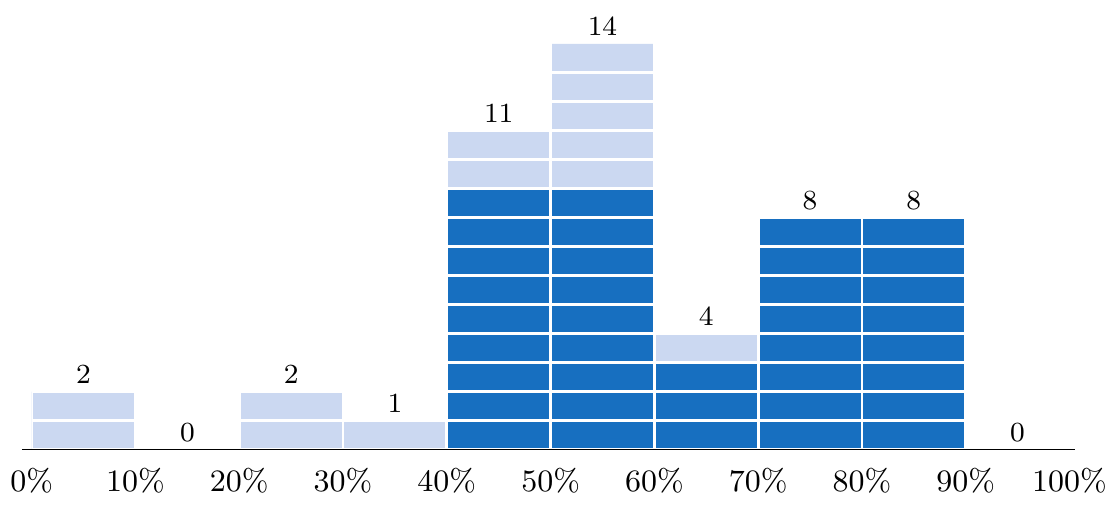}\hspace{10pt}\includegraphics[width=190pt]{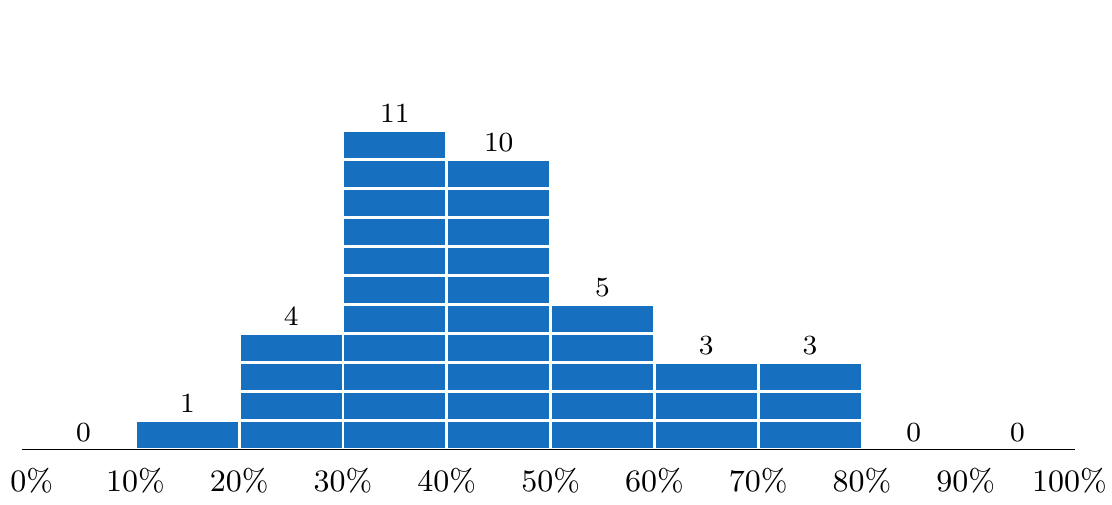}

\includegraphics[width=190pt]{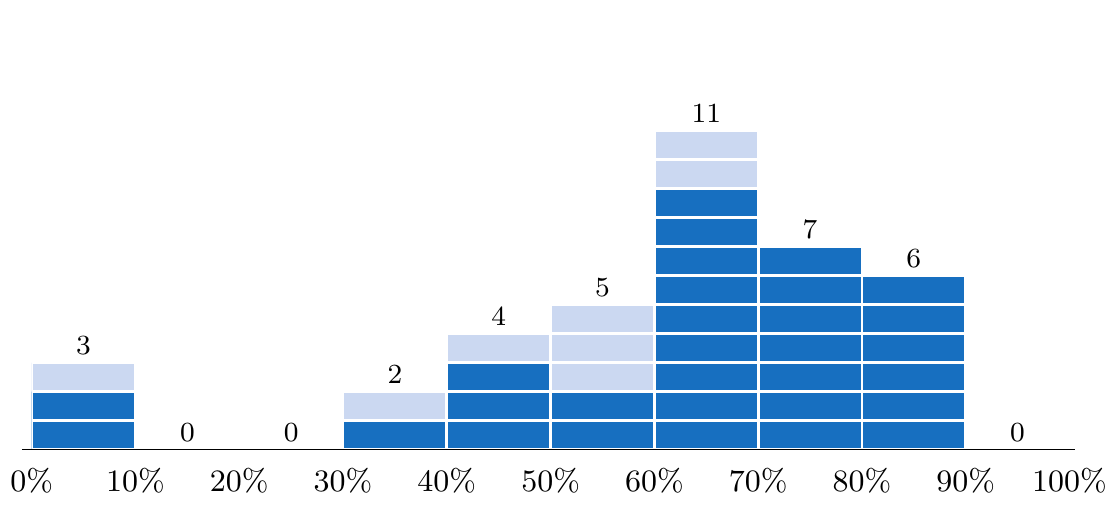}\hspace{10pt}\includegraphics[width=190pt]{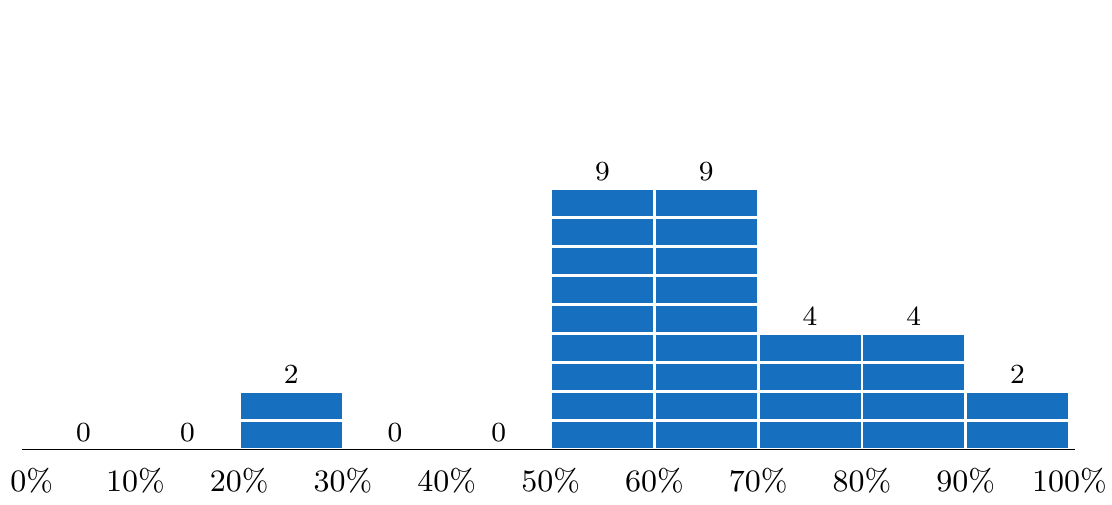}
\end{center}
\caption{Score distribution for homework (left) and exam (right) in the courses Mathematics 2 with teacher-graded homework (top) and Mathematics 3 with self-graded homework (bottom).  Light blue colored boxes stand for students that submitted homework but did not participate in the exam.}\label{FIG-Histo}
\end{figure}

In order to analyze our data further, I will now compare the distribution of homework and exam scores.  I start by comparing the percentage of tasks in Mathematics 3 with self-given label ``ok'' with the numerical exam score. This can be seen in the lower part of Figure \ref{FIG-Histo}. In comparison, the upper part of Figure \ref{FIG-Histo} shows the data from Mathematics 2, where the homework was teacher-graded. Light blue boxes represent students who submitted homework but were either not admitted to the exam or were admitted but did not participate for personal reasons.

\vspace{5pt}

\smallskip

While the teacher graded homework shows a U-shaped distribution with many students scoring just short above the 40\% threshold, the self-graded homework has a more Gaussian-like shape.\footnote{Indeed, the two dark blue boxes on the very left belong to two students, who submitted homework throughout the term but stopped submitting self-evaluations a few weeks in.} Although with the self-graded homework more students ``scored'' below the 40\%, the exam results indicate that relatively more students arrived finally at a pass, while in the semester with teacher-graded homework the exam distribution is Gaussian-like with more students failing the exam compared to Mathematics 3.
\medskip

Now I would like to survey gender aspects of my approach. In Mathematics 2 we had 25 male students and 25 female students that submitted homework in seven groups with two male students, eight groups of two female students, nine groups of one male and one female student. Two male students submitted alone. In Mathematics 3 we had 19 male students and 19 female students who submitted in six groups of two male students, six groups of two female students, five groups of one female and one male student. Moreover, we had one group with two female and one male students, and one male student who submitted alone.

\smallskip

Table \ref{TAB-3} shows the means and standard deviations of all/female/male students' scores in homework and exam. As above, the homework score is the score of the respective team. In Mathematics 3 it is the percentage of self-given label ``ok.'' As mentioned earlier, not all students who submitted homework (and self-evaluations) did participate in the exam. In order to compare for instance the mean score in the exam with the mean score in the homework,  Table \ref{TAB-3} shows the homework mean only of those students who participated in the exam. The values of all students that submitted homework is additionally indicated in square brackets. The same goes for the standard deviations.

\smallskip

\begin{table}
\begin{center}
\begin{tabular}{lcccc}
\toprule
&\multicolumn{2}{c}{Homework}   &\multicolumn{2}{c}{Exam} \\
      \cmidrule(lr){2-3}      \cmidrule(lr){4-5}
   &  \small \hspace{15pt}Mean\hspace{15pt} &\small \hspace{18pt}SD\hspace{18pt}  & \small \hspace{15pt}Mean\hspace{15pt} &\small \hspace{18pt}SD\hspace{18pt}  \\
   \cmidrule(lr){2-5} 
  & \multicolumn{4}{c}{\small Mathematics 2 (teacher-graded homework)} \vspace{3pt} \\
all ($n=37\,[50]$)    &63\%\,[57\%] 	&15\%\,[19\%]  & 45\%  & 15\% 	\\
male ($n=19\,[25]$)   &59\%\,[53\%] 	&16\%\,[21\%]  & 46\%  & 14\%	\\
female ($n=18\,[25]$) &68\%\,[62\%] 	&14\%\,[17\%]  & 44\%  & 16\% 	\\
\midrule
  & \multicolumn{4}{c}{\small Mathematics 3 (self-graded homework)}  \vspace{3pt}\\
all ($n=30\,[38]$)    &63\%\,[59\%] 	&20\%\,[21\%]  & 65\%       & 16\%\\
male ($n=15\,[19]$)   &66\%\,[60\%]	    &16\%\,[22\%]  & 69\%       & 17\% 	\\
female ($n=15\,[19]$) &59\%\,[58\%]     &24\%\,[21\%]  & 62\%       & 16\% 	\\
\bottomrule
\end{tabular}
\end{center}

\caption{Means and standard deviations of teacher-graded and self-graded homework as well as exam results. The numbers in brackets refer to mean and standard when students who did not participate in the exam are included.}\label{TAB-3}

\end{table}

Although there are differences between male and female students, they are very small and likely to be due to the fact that our sample is rather small. Notice that our data, in particular the Mathematics 2 homework, can not be regarded as Gaussian which means that it is impossible to discuss if differences are significant. One can see nevertheless that in both exams men scored slightly better than women. In the teacher-graded homework this was reversed but not in the self-graded homework. I emphasize here that one should be very cautious when interpreting the percentage of self-given labels ``ok'' as a score.

\smallskip

\begin{figure}
\begin{center}
\includegraphics[width=190pt]{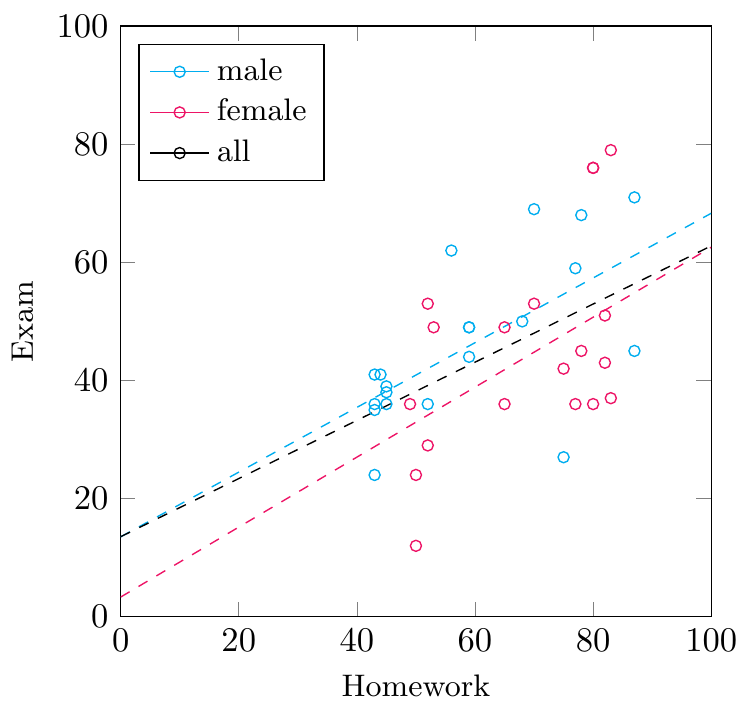}\hspace{5pt}\includegraphics[width=190pt]{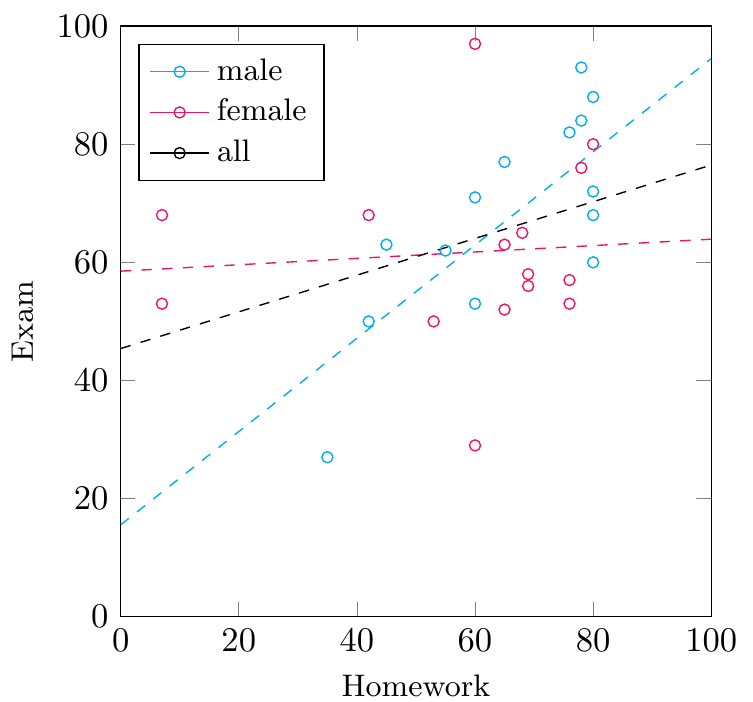}
\end{center}
\caption{Homework and exam scores in the teacher-graded course Mathematics 2 (left) and in the course Mathematics 3 with self-graded homework (right). Dashed lines represent linear regressions.}\label{FIG-Gender}
\end{figure}

In Figure \ref{FIG-Gender}, I plotted exam score over homework score for all students where I had both values (notice that this excludes students who did not get the exam admission and those who got admitted but did not participate in the exam). The dashed lines show linear regressions.  Notice that in the Mathematics 2 exam the threshold to pass was 35\% while in the Mathematics 3 exam it was 50\% and that the exams were designed by different teachers. The two red circles on the very left of the right picture belong to the two students mentioned earlier that submitted homework throughout the semester but stopped submitting self-evaluations a few weeks in.

\smallskip

Taking gender into account one can see that in the fully teacher-graded setting there is no qualitative difference to be detected. In contrast to that the setting with self-graded homework shows large differences between male and female students. Since the data on the right is much more ``spread out'' it is difficult to draw any general conclusions. One can see however that the best exam (97\%) was written by a female student who gave herself only ``ok'' on 60\% of the homework. This suggests that she used the label ``needs work'' to revise effectively\footnote{Indeed, this particular student gave herself several times ``needs work'' also for submissions that were quite close to perfect.}. On the other hand the right upper corner of the right picture suggests that among those students that scored between 60\% and 80\% in the homework the male ones achieved often higher scores in the exam. 

\medskip

I conclude this section with two feedback comments from the official and anonymous course evaluations that usually takes place halfway into the semester. The first comment was received in the course Mathematics 3, that is, before any exam was written, and during the very first experience with ungrading.

\vspace{5pt}

\begin{center}
\begin{minipage}{350pt}
\textquotedblleft{}The ungrading system motivated me very much. Deductions of points in homework usually made me angry/disappointed/more demotivated. Now, I engage much more with the teacher's feedback. This increases my learning progress.\textquotedblright{}
\end{minipage}
\end{center}

\vspace{5pt}

The second comment was received in the course Mathematics 4, i.e., after the completion of Mathematics 3 and the exam, and thus reflects on the whole ungrading experience from that course. Recall that in Mathematics a system without points on homework was maintained but teacher-given feedback of the form ``exemplary,'' ``correct,'' ``half-correct,'' ``wrong'' was given. Most importantly, there was no threshold for exam participation.

\vspace{5pt}

\begin{center}
\begin{minipage}{350pt}
\textquotedblleft{}After all, the 9 credits are more than well deserved. Compared with other courses, where we, for example, get 6 credits, we have to work for this one twice as much. That there is no more points on the homework reduces however the pressure. I consider that to be very positive.\textquotedblright{}
\end{minipage}
\end{center}

\vspace{5pt}

The latter two comments indicate that the ungrading approach reduced students' stress and fostered a focus on learning rather than scoring when doing homework assignments. This is contrasted by the light blue boxes in Figure \ref{FIG-Histo} where one can see that seven students scored above the exam admission threshold but did not attempt the exam. Since by study regulations every exam can only be attempted at most four times, this suggests that students did not write the exam because they saw only a small chance of passing.

\subsection{Reflecting on the Ungrading Experience by Dr.\  Wegner}

During the semester I experienced some of the typical effects of ungrading as mentioned in the literature: Less stressed students, no discussions/complains about points, more questions of the type ``How do I do that right?'' instead of ``Why didn't I get the points?'' I saw students after the lesson sitting together and discussing what they should write for their self-evaluation and indeed there were some good comments that explain the mathematical problem and how it was overcome during the recitation. On the other hand, I saw also that a majority of students did not make much of the approach: The percentage of insightful comments was rather small and as outlined in the case analysis in section \ref{SEC-Sven-2}, some tasks (e.g. the Taylor polynomial) triggered no insightful comment although it would have been very easy and could have been very beneficial for students to write something meaningful. This feels disappointing for me as their teacher.

\smallskip

During the semester, some students expressed concerns about the consultations since in contrast to earlier years there was no weekly update on their performance in terms of grades for each assignment. I managed however to convince them to trust their own judgement---and to trust me to play fair in the consultation. Doing the consultations felt very good and my impression was that many students benefited from this (see the uniform continuity task in section \ref{SEC-Sven-2}). Finding time for $17\times 40$ minutes in the last week of the term was however difficult and caused a much bigger workload than the traditional setup with maybe one hour looking at scores and one or two critical cases and then a two hours drop-in session before the exam for questions. 

\smallskip

Due to the constraints in the study regulations (for instance that courses need to be completed by an exam) or rule-like customs (for example that it is very common that students submit weekly homework and are only allowed to write the exam if they score above a certain threshold) it was not possible for me to implement ungrading in a way described by Dr.\ von Renesse in section \ref{SEC-Ungrading-WSU}. On the other hand, most German universities do not ``grade behavior'' anyway, e.g., it is not common to take attendance into account. The fact that there are usually no weekly quizzes or midterm exams, is also in favor of students who start with lower premises but improve throughout the course. Nevertheless I think that the ungrading elements that I used, helped students to grow confidence, that it reduced stress and that it helped students to move their focus from scoring to understanding. The latter is an emancipation process that all students should go through anyway (and almost all do), but my feeling is, that ``limited ungrading'' in the form that I used can support and possibly accelerate this process which is beneficial for students and for teachers. In view of the large amount of extra work I am however not sure if I would go for this model again.

\smallskip

Concerning the data that I collected and analyzed, many limitations have already been mentioned: My sample was rather small and when comparing two consecutive courses with different topics it is very difficult to draw reliable conclusions. Comparing my students' self-evaluations with teacher-given grades is questionable since the premise of the self-evaluation was not to emulate the teacher's grading process but to bookmark tasks that need further attention before the exam. Concerning the gender-aspects I can only say that it is likely that there are differences between the approaches of male and female students, but it would need a much more sophisticated study to understand the details.

\section{CONCLUSION} 

The German and the U.S.\ education systems are very different at the tertiary level, in terms of student preparation, students' expectations and general set-up. In addition, we teach at very different institutions in terms of size and faculty expectations. Therefore the possibilities to use ungrading are different as well. Interestingly, some of the results are similar though: there are clear benefits like students' motivation and decreased stress levels, drawbacks like additional workload for professors, as well as puzzling patterns in the gender analysis. In Dr.\  von Renesse's classes, female students do overall more revisions, but advocate less for themselves during final conferences and grade discussions. In Dr.\  Wegners' class, no clear correlation between homework and exam results be detected when we consider gender as a separate variable. 

In the literature about inquiry-oriented learning, we can see that female students could be in a disadvantaged situation compared to a traditional lecture system, see \cite{Lubienski, Reinholz}. Both teaching with inquiry and ungrading increase the possibilities of including bias, both from the students and from the professors. A larger study is needed to examine these patterns more. Being aware should help us though to inform our teaching choices. Earlier journal assignment can, for instance, open the conversation about gender difference in portrayed confidence and self advocacy. We also need to consider that students are used to being graded and being told what to do next. They need to practice reflection, self-evaluation, and how to make independent, informed decisions to improve. Our assignments need to support students in their growth toward these meta goals.

\providecommand{\href}[2]{#2}

\end{document}